\documentclass[]{interact}

\usepackage{epstopdf}
\usepackage[caption=false]{subfig}
\usepackage{pgfplots}
\pgfplotsset{compat=1.16}

\usepackage[numbers,sort&compress]{natbib}
\bibpunct[, ]{[}{]}{,}{n}{,}{,}
\makeatletter
\def\NAT@def@citea{\def\@citea{\NAT@separator}}
\makeatother

\theoremstyle{plain}
\newtheorem{theorem}{Theorem}[section]
\newtheorem{lemma}[theorem]{Lemma}
\newtheorem{corollary}[theorem]{Corollary}
\newtheorem{proposition}[theorem]{Proposition}
\newtheorem{conjecture}[theorem]{Conjecture} 

\theoremstyle{definition}
\newtheorem{definition}[theorem]{Definition}
\newtheorem{example}[theorem]{Example}

\theoremstyle{remark}

\usepackage{hyperref}
\usepackage{xcolor}
\begin{document}
\title{Box dimension of stable sub-slices of fractal graphs over Anosov automorphisms}

\author{
\name{Bernardo Carvalho\textsuperscript{a}\textsuperscript{b}\thanks{CONTACT Bernardo Carvalho. Email: bmcarvalho@lncc.br} and Rafael da Costa Pereira\textsuperscript{c}\thanks{CONTACT Rafael da Costa Pereira. Email: rafaeldacostapereira@gmail.com}}
\affil{\textsuperscript{a} National Laboratory of Scientific Computing, Av. Getúlio Vargas 333, CEP 25651-070, Petrópolis – RJ, Brazil; \textsuperscript{b} Dipartimento di Matematica, Universit\`a degli Studi di Roma Tor Vergata, Via Cracovia n.50 - 00133, Roma - RM, Italy; \textsuperscript{c}Departamento de Matem\'atica, Universidade Federal de Minas Gerais - UFMG, Av. Ant\^onio Carlos, 6627 - Campus Pampulha, 
Belo Horizonte - MG, Brazil}
}

\maketitle

\begin{abstract}
We consider fractal graphs invariant by a skew product $F\colon\mathbb{T}^k\times \mathbb{R}\rightarrow \mathbb{T}^k\times \mathbb{R}$ of the form $F(x,y)=(Ax, \lambda y+p(x))$ where $0<\lambda<1$, $p\colon\mathbb{T}^k\to\mathbb{R}$ is a $C^{k+1}$ function, and $A$ is an Anosov automorphism of $\mathbb{T}^k$ admitting $k$ distinct eigenvalues with respective eigenvectors forming a basis of $\mathbb{R}^k$. We extend the techniques introduced by \citet*{kaplan1984} to calculate the box dimension of \emph{sub-slices} of the invariant graph of $F$. We show that the stable sub-slices have a different fractal structure in comparison to the graph. We exhibit conditions on the skew product that ensure the box dimension of the graph is smaller than the sum of the box dimensions of its stable/unstable sub-slices. These conditions hold for generic functions $p\in C^{k+1}$.
\end{abstract}

\section{Introduction}

\hspace{+0.5cm}Fractal sets appear naturally in hyperbolic dynamics, as typically, hyperbolic sets exhibit interesting fractal structures. The study of these properties, utilizing metrics such as box or Hausdorff dimension, was particularly influenced by the classical example of the Weierstrass function, which served as a prominent source of inspiration for early dimensional calculations \cite{falconer2014fractal}. In the context of hyperbolic dynamics, analyzing the stable and unstable manifolds offers a pathway to calculate the fractal properties of hyperbolic basic sets. For instance, in surface dynamics, the dimension of a hyperbolic set can often be understood through the dimensions of its local stable and unstable slices — the intersections of the local stable/unstable manifolds with the set itself \cite{manning1983,Takens,Palis-Viana}.

Extending this result to higher dimensional (and non-conformal) hyperbolic sets is an open problem that may require the development of new tools. The lack of regularity on the stable/unstable holonomies and the possible difference between Hausdorff and box dimensions are difficulties that appear in this scenario. 

In this direction, we remark the work \cite{hasselblatt2004} in which the authors prove that the ideas of \cite{manning1983} in the case of surface diffeomorphisms can be applied to the Smale Solenoid, which is a hyperbolic set on a manifold of dimension three, and state the following conjecture:
\begin{conjecture}
The fractal dimension of a hyperbolic set is (at least generically or under mild hypotheses) the sum of those of its stable and unstable slices, where ``fractal'' can mean either Hausdorff or upper box dimension.
\end{conjecture}

While our primary focus is on calculating the box dimension, we note that it serves as a natural upper bound for the Hausdorff dimension, and the two values coincide in some cases, such as in the case of the graph of the Weierstrass function. Recently, the Hausdorff dimension of the graph of the classical Weierstrass function was calculated in \cite{BARANSKI201432}, using techniques from \cite{L-YOUNG} and \cite{Tsujii_2001}. A different proof of the same result using a more direct approach via conditional probabilities was obtained in \cite{Keller}. This was also generalized by \cite{Shen}, replacing the cosine in the Weierstrass function with a $C^2$-periodic function. 

Some recent advances in the above-mentioned conjecture were obtained in \cite{MOHAMMADPOUR_PRZYTYCKI_RAMS_2022} for three-dimensional thin solenoids — those whose contractions of fibers of the skew product defining it are so strong that stable slices have dimension strictly less than one — while \cite{Bortolotti_2022} extended the result to higher-dimensional thin solenoids. Additionally, the conjecture was established for the particular case of $C^1$ average conformal hyperbolic sets \cite{JuanWang}.

Another significant contribution to the theory was obtained by Kaplan, Mallet-Paret, and Yorke in \cite{kaplan1984} where the box dimension of fractal graphs invariant by a skew product over the cat map on $\mathbb{T}^2$ was computed using Fourier analysis of almost periodic functions. These techniques are quite distinct and non-standard and will be explained in detail in Section \ref{T2} of this article. These graphs can be defined following \cite{kaplan1984} considering skew products of the form $F\colon\mathbb{T}^2\times \mathbb{R}\rightarrow \mathbb{T}^2\times \mathbb{R}$ defined by
\begin{equation}\label{F}
F(x,y)=(Ax, \lambda y+p(x)),
\end{equation}
where $A$ is an Anosov automorphism, $0<\lambda<1$, and $p\colon\mathbb{T}^2\to\mathbb{R}$ is a smooth function. This skew-product has an attractor that can be seen as a graph over $\mathbb{T}^2$ (see \cite{kaplan1984}) such that any point $(x,y)$ in the attractor satisfies
\[
y=\sum_{i=1}^{\infty}\lambda^{i-1}p(A^{-i}x).
\]
For simplicity, the authors in \cite{kaplan1984} consider the function $\phi\colon\mathbb{T}^2\to\mathbb{R}$ defined by
\begin{equation}
\phi(x)=\sum_{n=0}^{\infty} \lambda^n p(A^{-n}x).
\label{phi}
\end{equation}
Thus, $y(x)=\phi(A^{-1}x)$, and since $A$ is an Anosov automorphism, the box dimensions of the graphs of $y$ and $\phi$ are the same. The calculated box dimension of the stable slice of the graph of $\phi$ equals the box (and Hausdorff) dimension of the graph of the Weierstrass function. 

It is noted in \cite{kaplan1984} that their techniques of Fourier analysis of almost periodic functions can be used to calculate the box dimension of fractal graphs over higher-dimensional Anosov automorphisms of $\mathbb{T}^k$ (see Theorem C in \cite{kaplan1984}). The argument is similar to the case of $\mathbb{T}^2$ using the strongest stable eigenvalue of the base hyperbolic matrix to control the sizes of boxes that cover the graph. As is expected, see for instance the works \cite{Stark,HADJILOUCAS_NICOL_WALKDEN_2002,Katrin}, the regularity of the graph depends only on the relation of the contraction on the fiber of the skew product and the strong stable direction of the base automorphism, that is, the values of the other stable eigenvalues and also their respective eigenvectors do not influence the box dimension of the graph. 

The primary contribution of this article is to place the constructions introduced in \cite{kaplan1984} into a broader, higher-dimensional context. We extend the techniques introduced there to calculate the box dimension of \emph{sub-slices} of invariant graphs over higher-dimensional Anosov automorphisms of $\mathbb{T}^k$. For the skew product $F\colon\mathbb{T}^k\times \mathbb{R}\rightarrow \mathbb{T}^k\times \mathbb{R}$ defined as in (\ref{F}), we assume that $A$ has $k$ distinct eigenvalues $B_1,\dots,B_k$ and that the respective eigenvectors $v_1,\dots,v_k$ form a basis of $\mathbb{R}^k$. Write points of $\mathbb{T}^k$ on this basis as $\underline{t}=(t_1,\ldots,t_k)$ and write $\phi$ on this basis as
\[
\phi(\underline{t})=\sum_{n=0}^{+\infty} \lambda^n p(t_1B_1^{-n}v_1+\cdots+t_kB_k^{-n}v_k)=\sum_{n=0}^{+\infty} \lambda^n p(t_1B_1^{-n},\dots,t_kB_k^{-n}).
\]

\begin{figure}[!ht]
    \centering
    \includegraphics[width=0.5\linewidth]{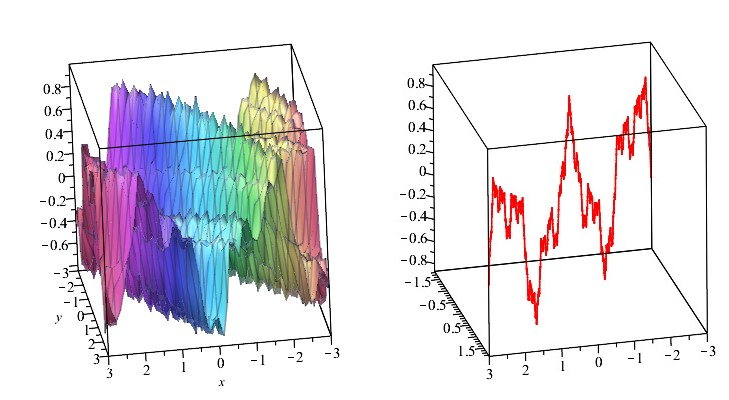}
    \caption{The fractal graph (left) and the graph of a sub-slice (right).}
    \label{slice e subslice}
\end{figure}

At any point of the graph $(\underline{t},\phi(\underline{t}))$, each direction $v_i$ gives a slice of the graph given by the graph of the function $\phi_i$ defined by
\[
\phi_i(t):=\phi(tv_i+\underline{t}),
\]
with $t$ sufficiently small. The maps $(\phi_i)_{i=1}^k$ are called the \emph{sub-slices} of the graph $\phi$ at the point $(\underline{t},\phi(\underline{t}))$ (see Figure \ref{slice e subslice}).

By decomposing the stable and unstable slices into these sub-slices — each corresponding to a distinct eigenspace of the base automorphism — we establish conditions under which the box dimension of the full graph is strictly less than the sum of the box dimensions of its stable and unstable sub-slices. This constitutes our main result (Theorem 3.6), demonstrating that the fractal structure of the sub-slices can differ fundamentally from the overall structure of the graph. We analyze the genericity of these conditions (see Theorem \ref{3.9}) and present an example of a graph invariant by a skew product satisfying such conditions in Example \ref{example}.

We emphasize that our results apply to Anosov automorphisms on the base. This restriction is necessary because the methods we adapt from \cite{kaplan1984} rely on having a uniform, global eigenspace structure. Since the matrix $A$ and its eigenvectors are constant, we can separate the dynamics along distinct linear directions and apply Fourier analysis. If we considered non-linear Anosov diffeomorphisms, the invariant directions would vary from point to point, preventing the clear separation needed to compute the dimensions in this way. These observations raise open problems for future investigation. Extending our dimensional analysis to non-linear Anosov diffeomorphisms would likely require different techniques to properly handle non-constant invariant bundles. Furthermore, even in the linear case, the calculation of the box dimension for sub-slices along non-invariant directions remains open, as the lack of dynamical decoupling introduces additional geometric complexities.

\section{Kaplan, Mallet-Paret, and Yorke techniques}

In this section, we define the basic concepts of fractal geometry and hyperbolic dynamics that will be used in this article and explain the techniques of \cite*{kaplan1984} to calculate the box dimension of the graph invariant by $F$.

\subsection{Box dimension}
In this subsection, we present some basic results about the Box dimension.
\begin{definition}
Consider $X$ a non-empty totally bounded subset of a metric space $M$ and for each $\delta>0$ 
let $N_{\delta}(X)$ be the smallest number of sets of diameter at most $\delta$ which union covers X. We define the lower box dimension of $X$ as
\begin{align*}
\text{\underline{dim}}(X)=\liminf_{\delta \rightarrow 0^+} \frac{\log{N_{\delta}(X)}}{-\log{\delta}}
\end{align*}
and the upper box dimension of $X$ as
\begin{align*}
\overline{\text{dim}}(X)=\limsup_{\delta \rightarrow 0^+} \frac{\log{N_{\delta}(X)}}{-\log{\delta}}.
\end{align*}
If they are equal, this is the box dimension of $A$, and we write $\text{dim}(X)$.
\end{definition}
The following proposition exhibits an easy way to calculate the box dimension of a set and will be used a few times in the arguments. This is a classical fact of the box dimension and we include its proof for completeness.

\begin{lemma}\label{box}
If there are positive constants $C_1,C_2$ and $a$ such that
\[
C_1\delta^{-a}\leq N_{\delta}(X) \leq C_2\delta^{-a},
\]
then $\text{dim}(X)=a$.
\end{lemma}

\begin{proof}
We just need to note the following inequalities:
\begin{eqnarray*}
C_1\delta^{-a}\leq  N_{\delta}(X) \leq C_2\delta^{-a} &\Rightarrow & \frac{\log (C_1\delta^{-a})}{-\log \delta}\leq \frac{\log N_{\delta}(X)}{-\log\delta} \leq \frac{\log(C_2\delta^{-a})}{-\log \delta} \\
&\Rightarrow & \frac{\log C_1}{-\log \delta}+a\leq \frac{\log N_{\delta}(X)}{-\log\delta} \leq  \frac{\log C_2}{-\log \delta}+a.
\end{eqnarray*}
Then $\text{dim}(X)=a$ follows by letting $\delta\to0$ in the last inequalities.
\end{proof}

In the special case that $X$ is the graph of a continuous function $f\colon I=[c,d]\to\mathbb{R}$, Lemma \ref{box} can be applied to calculate the box dimension considering similar inequalities for the variation of the function $f$ defined by 
$$\underset{I}{\text{var}}(f)=\sup_{t \in I} f(t) - \inf_{t \in I} f(t).$$ 
The following results can be obtained as particular cases of Lemma 7 from \cite{Przytycki1989}. We include a proof for completeness.

\begin{proposition}\label{var}
If there are $L_0\in(0,1)$, $a\in(0,1)$, $C_1,C_2>0$ such that
\[
C_1 L^a \leq \underset{J}{\text{\upshape var}}(f) \leq C_2 L^a
\]
for every interval $J$ with length$(J)=L\leq L_0$, then
\[
\text{\upshape dim(graph}(f))= 2-a.
\]
\end{proposition}
\begin{proof}
Given $\delta<\text{min}\{L_0,d-c,1\}$, let 
\[
|I|_\delta:= \left\lceil \frac{d-c}{\delta}\right\rceil,
\]
where $\lceil x\rceil$ is the smaller integer greater than $x$.
For each $j\in\{1,\dots, |I|_\delta \}$, let
$$I_j:=[c+(j-1)\delta,c+j\delta], \,\,\,\,\,\, |f|_j:= \left\lceil \frac{\underset{I_J}{\text{var}}(f)}{\delta}\right\rceil,$$
and for each $k\in \{1,\dots, |f_j|\}$, let 
$$I^k_{j}=\left[\inf_{t \in I_j} f(t)+(k-1)\delta, \inf_{t \in I_j} f(t)+k\delta\right].$$ Consider the set 
$$I=\{I_j\times I^k_{j}; j\in\{1,\dots, |I|_\delta\}, k\in \{1,\dots, |f_j|\}\}$$ and let 
\[
A:=\#I=\sum_{k=1}^{|I|_\delta} |f_j|.
\]
The set $I$ covers $\text{graph}(f)$ with sets of diameter $\sqrt{2}\delta$, so 
\[
N_{\sqrt{2}\delta}(\text{graph}(f))\leq A. 
\]

\begin{figure}
    \centering
{\resizebox*{5cm}{!}{\includegraphics{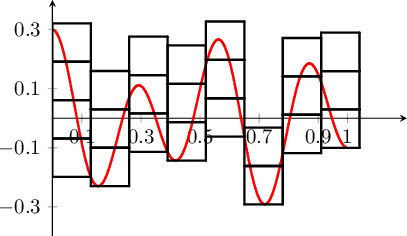}}}
    \caption{Example of the cover $I$ for $\delta=0.13$ and $f(x)=0.2 \cos(7 \pi x)+0.1 \cos(4 \pi x)$ on $[0,1]$.}
    \label{fig:enter-label}
\end{figure}

Note that $|I|_\delta\leq \frac{d-c}{\delta}+1$ and that by hypothesis
\[
\underset{I_J}{\text{var}}(f)\leq C_2\delta^a \,\,\,\,\,\, \text{for every} \,\,\,\,\,\, j\in\{1,\dots, |I|_\delta\}.
\]
It follows that
\[
|f|_j\leq C_2\delta^{a-1}+1 \,\,\,\,\,\, \text{for every} \,\,\,\,\,\, j\in\{1,\dots, |I|_\delta\}
\]
and, hence,
\begin{align*}
N_{\sqrt{2}\delta}(\text{graph}(f))\leq A&\leq \left(\frac{d-c}{\delta}+1\right)\left(C_2\delta^{a-1}+1\right)\\ &=\delta^{a-2}(C_2(d-c)+(d-c)\delta^{1-a}+C_2\delta+\delta^{2-a}) \\ &\leq \delta^{a-2}(C_2(d-c)+(d-c)+C_2 +1) \\ &\leq (\sqrt{2})^{a-2}\delta^{a-2}K_2
\end{align*}
for some positive constant $K_2$.
Also note that 
\[
N_{\sqrt{2}\delta}(\text{graph}(f))\geq \frac{A}{9}
\]
since any set of diameter $\sqrt{2}\delta$ intercepts no more than 9 distinct sets of the form $I_j\times I^j_{k}$. By hypothesis we have
\[
\underset{I_J}{\text{var}}(f)\geq C_1\delta^a \,\,\,\,\,\, \text{for every} \,\,\,\,\,\, j\in\{1,\dots, |I|_\delta\}
\]
and since we similarly have $|I|_\delta\geq \frac{d-c}{\delta}$ and
\[
|f|_j\geq C_1\delta^{a-1} \,\,\,\,\,\, \text{for every} \,\,\,\,\,\, j\in\{1,\dots, |I|_\delta\},
\]
it follows that
\[
N_{\sqrt{2}\delta}(\text{graph}(f))\geq\frac{A}{9}\geq \left(\frac{d-c}{9\delta}\right)C_1\delta^{a-1}\geq (\sqrt{2})^{a-2}\delta^{a-2}K_1
\]
for some positive constant $K_1$. Thus,
\[
K_1(\sqrt{2}\delta)^{a-2}\leq N_{\sqrt{2}\delta}(X) \leq K_2(\sqrt{2}\delta)^{a-2}
\]
and Proposition \ref{box} ensures that $\text{dim}(\text{graph}(f))=2-a$.
\end{proof}

The essence of this proof is that the number of intervals of size $\delta$ needed to cover the domain are proportional to $\delta^{-1}$ and on each of these intervals the number of intervals of size $\delta$ needed to cover the image is, by hypothesis, proportional to $\delta^{a-1}$, so the number of squares of size $\delta$ needed to cover the graph is proportional to $\delta^{a-2}$. This can also be done for a $n$-variable function with the only difference being that the number of $n$-cubes needed to cover the domain is proportional to $\delta^{-n}$. Then we can state the following result.

\begin{proposition}\label{nvar}
If $f\colon I^n\to\mathbb{R}$ is a continuous function and there are $L_0\in(0,1)$, $a\in(0,1)$, $C_1,C_2>0$ such that
\[
C_1 L^a \leq \underset{J}{\text{\upshape var}}(f) \leq C_2 L^a
\]
for every $n$-cube $J$ with side $L\leq L_0$, then
\[
\text{\upshape dim(graph}(f))= n+1-a.
\]
\label{grafico}
\end{proposition}

\subsection{Hyperbolic Sets and Stable/Unstable Manifolds}

The following is the classical definition of a hyperbolic set.
\begin{definition}[Hyperbolic set] Let $M$ be a Riemannian manifold, $U \subset M$ be an open set, and $f\colon U\rightarrow M$ be a smooth embedding. A compact invariant subset $\Lambda \subset U$ is called a \textit{hyperbolic set} if the tangent bundle over $\Lambda$ splits as $T_{\Lambda}M=E^s\oplus E^u$ and there exist constants $C>0$ and $0<\lambda<1$ such that $\|Df_{\vert_{E^s}}^n\| <C\lambda^n$ and $\|Df_{\vert_{E^u}}^{-n}\| <C\lambda^n$ for every $n \in \mathbb{N}$.
\end{definition}
The Stable Manifold Theorem ensures that if $\epsilon$ is small enough, then the sets
\begin{align*}
W^s_{\epsilon}(x)&=\{y \in M | \text{ dist}(f^{n}(x),f^{n}(y))\leq\epsilon   \text{ for all }   n \in \mathbb{N}\} \\
W^u_{\epsilon}(x)&=\{y \in M | \text{ dist}(f^{-n}(x),f^{-n}(y))\leq\epsilon   \text{ for all }   n \in \mathbb{N}\}
\end{align*}
are $C^1$ embedded disks tangent to $E^s$ and $E^u$, respectively. In particular, there is a local product structure close to any $x\in\Lambda$: there exists $\delta>0$ and a map $$[,]\colon (W^s_{\epsilon,\delta}(x)\cap\Lambda)\times (W^u_{\epsilon,\delta}(x)\cap\Lambda)\to M,$$ where $W^s_{\epsilon,\delta}(x)=W^s_{\epsilon}(x)\cap B(x,\delta)$ and $W^u_{\epsilon,\delta}(x)=W^u_{\epsilon}(x)\cap B(x,\delta)$, defined by $$[p,q]=W^{u}_{\epsilon}(p)\cap W^{s}_{\epsilon}(q).$$
Thus, for each $x\in\Lambda$, the map $[,]$ defines a homeomorphism between the product $$(W^s_{\epsilon,\delta}(x)\cap\Lambda)\times (W^u_{\epsilon,\delta}(x)\cap\Lambda)$$ and $\Lambda\cap B(x,\delta')$ (for some $\delta'\in(0,\delta]$).
The sets in this product are called the local stable and local unstable slices of $\Lambda$ at $x\in\Lambda$.
If the map $[,]$ is bi-Lipschitz, then the dimension of $\Lambda$ equals the sum of the dimensions of the stable and the unstable slices at any point, since Hausdorff and box dimensions are preserved by bi-Lipschitz functions. This fact is used in \cite{manning1983} to calculate the dimension of two dimensional horseshoes.
\begin{theorem}[Theorem 2 of \citep{manning1983}]
Let $f\colon S\rightarrow S$ be a $C^2$ axiom A diffeomorphism of a surface, and $\Lambda\subset S$ be a basic set of $f$. Then the dimension of $\Lambda$ is the sum of the dimensions of the stable and unstable slices of any $x\in\Lambda$. Also, the dimension varies continuously with $f$.
\end{theorem}

Here dimension means both box and Hausdorff since these are equal in the case of hyperbolic basic sets of surface diffeomorphisms. For hyperbolic sets on manifolds of higher dimensions, the map $[,]$ is always Hölder continuous but not necessarily Lipschitz (see \cite{shub2013global}). Thus, new techniques must be necessary to compute the dimension in this case.



\vspace{+0.4cm}

\subsection{Regularity of the graph and its sub-slices}

The global attractor of the skew product $F$ is represented by the graph of the function $\phi \colon \mathbb{T}^k \to \mathbb{R}$, which forms a $k$-dimensional continuous surface embedded in $\mathbb{T}^k \times \mathbb{R}$. Since the base dynamics is governed by the Anosov automorphism $A$, we can analyze the fractal properties of this overall surface by decoupling it along its invariant directions. For a fixed base point $\underline{t}\in\mathbb{T}^k$ and a chosen eigenvector $v_i$ associated with the eigenvalue $B_i$, the sub-slice $\phi_i(t) = \phi(tv_i + \underline{t})$ effectively isolates a one-dimensional curve on this surface.

With the sub-slices of the graph formally defined, we now turn to determining their respective box dimensions. The regularity of each sub-slice $\phi_i$ depends on the competition between the contraction rate on the fiber, $\lambda$, and the dynamical behavior of the base automorphism along the corresponding eigenspace, given by $|B_i|$. We first address the regime where $\lambda < |B_i|$. Dynamically, this condition implies that the contraction along the fiber is strictly stronger than the contraction (or expansion) along the base direction $v_i$. As a result, the rapid decay of the $\lambda^n$ terms in the defining sum suppresses the emergence of rough, high-frequency oscillations. The sub-slice thus inherits the smoothness of the perturbation function $p$, preventing the formation of a fractal structure. The following proposition makes this precise by showing that $\phi$ is continuously differentiable along such directions, which immediately implies that the box dimension of these sub-slices is exactly one.

\begin{proposition}\label{maior}
If $\lambda <|B_i|$, then $\phi$ is a $C^1$ function of $t_i$ and $|\partial_i \phi|$ is bounded independent of $(t_1,\dots,t_k)$.
\label{dim1}
\end{proposition}
\begin{proof}
The following function is the candidate for the derivative of $\phi$ on $(t_1,\dots,t_k)$ in the direction $v_i$:
\[
h(\underline{t})=\sum_{n=0}^{+\infty}\left(\frac{\lambda}{B_i}\right)^n\frac{\partial p}{\partial t_i}(t_1B_1^{-n},\dots ,t_kB_k^{-n})
\]
Let us prove it is, indeed, the case. For each $m\in\mathbb{N}$, let
\[
h_m(\underline{t})=\sum_{n=0}^{m}\left(\frac{\lambda}{B_i}\right)^n\frac{\partial p}{\partial t_i}(t_1B_1^{-n},\dots ,t_kB_k^{-n}).
\]
Note that $h_m$ converges uniformly to $h$ because $\frac{\partial p}{\partial t_i}$ is bounded and $\lambda<|B_i|$. Therefore
\[
\int_{0}^{t_i} h_m(\underline{t}) dt_i \underset{m\to\infty}{\longrightarrow} \int_{0}^{t_i} h(\underline{t}) dt_i.
\]
The notation $\int_{0}^{t_i} g(\underline{t}) dt_i$ means that we are integrating $g$ as a function of one variable where $\underline{t}$ varies only on the coordinate $t_i$. We hope this does not cause confusion with the number $t_i$ that is fixed on the point $(t_1,\dots,t_k)$ we consider at the beginning. On the other hand we have
\[
\int_{0}^{t_i} h_m(\underline{t}) dt_i=\sum_{n=0}^{m} \lambda^np(t_1B_1^{-n},\dots,t_k B_k^{-n}) \underset{m\to\infty}{\longrightarrow} \phi(\underline{t}),
\]
because $p$ is bounded and $\lambda<1$ it follows that 
\[
\int_{0}^{t_i} h_m(\underline{t}) dt_i\underset{m\to\infty}{\longrightarrow} \phi(\underline{t}),
\]
So we have that $\phi(\underline{t})=\int_{0}^{t_i} h(\underline{t}) dt_i$, and, hence,
\[
h(\underline{t})=\frac{\partial \phi}{\partial t_i}(t_1,\ldots,t_k)
\]
completing the proof.
\end{proof}
We end this subsection noting that this applies to all unstable sub-slices since, in this case, the eigenvalues have modulus bigger than 1, while $\lambda<1$.

\vspace{+0.4cm}

\subsection{\texorpdfstring{Box dimension for the stable slice in $\mathbb{T}^2\times\mathbb{R}$}{Box dimension for the stable slice in T2 X R}}
\label{T2}
In this section, we consider the case $k=2$, that is, the skew product is of the form $$F\colon\mathbb{T}^2\times\mathbb{R}\rightarrow \mathbb{T}^2\times\mathbb{R}.$$ In this case, at any point of the graph there is only one stable slice and we discuss how to calculate its box dimension following \cite{kaplan1984}. This will be used in Section 3 to calculate the box dimension of all stable sub-slices in the higher-dimensional case associated to directions $v_i$ with $|B_i|<\lambda$. Using the notation of the previous subsection, when $k=2$ the hyperbolic matrix $A$ has two eigenvalues $B_1<1$ and $B_2>1$. If $v_1$ is the eigenvector associated to $B_1$, then the stable slice at any point $\underline{t}\in\mathbb{T}^k$ is given by the graph of the function $\phi_1$ defined by 
\[
\phi_1(t):=\phi(tv_1+\underline{t}).
\]
Letting $q(t)=p(tv_1+\underline{t})$, the stable slice is given by the graph of the function
\[
f(t)=\sum_{n=0}^{\infty}\lambda^n q((B_1^{-1})^{n} t).
\]
The following theorem calculates the box dimension of this stable slice:
\begin{theorem}[Theorem A in \cite{kaplan1984}]
\label{kaplanthmA}
Let
\[
f(t)=\sum_{n=0}^{\infty}\lambda^n q(B^{n} t)
\]
where $0<\lambda<1$ and $B>1/\lambda$ and assume that
\begin{equation}
q(t)=\sum_{i=1}^N q_i\cos{(a_it + \theta_i)},
\label{functionq}
\end{equation}
for some real numbers $q_i,a_i,\theta_i$. Then either $f$ is $C^1$, and hence $\mathrm{dim(graph}(f))=1$, or
\[
\mathrm{dim(graph}(f)) = 2-\left|\frac{\ln{\lambda}}{\ln{B}}\right|.
\]
\end{theorem}
Actually, the conclusion holds for every function $q$ that is almost periodic and satisfies hypothesis (H1) (see Definitions \ref{qp} and \ref{H1}). In this subsection, we explain the proof of this statement. First, we recall the necessary definitions and explain their consequences.

\begin{definition}\label{qp}
A function $q:\mathbb{R}\rightarrow\mathbb{R}$ is called \emph{almost periodic} if for every $\epsilon>0$ there is $l(\epsilon)>0$ such that any real interval of length $l(\epsilon)$ contains a number $\tau$ satisfying
\[
|q(x)-q(x+\tau)|<\epsilon \,\,\,\,\,\, \text{for every} \,\,\,\,\,\, x\in \mathbb{R}.
\]
\end{definition}
Even though $q$ is not periodic it is possible to define a formal Fourier series (see \cite*{corduneanu})  such that
\[
q(t) \sim \sum_{a} q_ae^{iat}
\]
where
\[
q_a=\lim_{T\rightarrow \infty} \frac{1}{2T}\int_{-T}^{+T} q(t)e^{-iat}\hspace{0.5em} dt .
\]
There is only a countable set of $a\in\mathbb{R}$ such that $q_a\neq 0$.\footnote{Notice that if $q$ was periodic this definition would coincide with the standard definition of Fourier Series.} 
The main idea to prove Theorem \ref{kaplanthmA} is to formally consider
\begin{equation}\label{g}
g(t):=\sum_{n=-\infty}^{\infty} \lambda^{n}q(B^{n} t),
\end{equation}
that contains the function $f$ but also the negative part of the series.
Ignoring the fact that the function $g$ may not be well defined, if we try to calculate its Fourier series coefficients we get:
\begin{align*}
g_\sigma & =\lim_{T\rightarrow \infty} \frac{1}{2T}\int_{-T}^{+T}\sum_{n=-\infty}^{\infty} \lambda^{n}q(B^{n} t)e^{-i\sigma t}\hspace{0.5em} dt \\
&= \lim_{T\rightarrow \infty} \sum_{n=-\infty}^{\infty} \frac{\lambda^{n}}{2T}\int_{-T}^{+T}q(B^{n} t)e^{-i\sigma t}\hspace{0.5em} dt  \\
&= \lim_{T\rightarrow \infty} \sum_{n=-\infty}^{\infty} \frac{\lambda^{n}}{2T}\int_{-T B^{n}}^{+T B^{n}}\frac{q(s)e^{-i\sigma s B^{-n}}}{B^n} \hspace{0.5em} ds\\
&  =  \sum_{n=-\infty}^{\infty} \lambda^{n} \lim_{TB^n \rightarrow \infty}  \frac{1}{2TB^n}\int_{-T B^{n}}^{+T B^{n}}q(s)e^{-i\sigma B^{-n} s }\hspace{0.5em} ds \\
& =\sum_{n=-\infty}^{\infty} \lambda^{n}q_{\sigma B^{-n}}.
\end{align*}
Regardless of the lack of rigor in the calculation above, all that is truly necessary is that
\begin{equation}\label{gsigma}
g_\sigma:=\sum_{n=-\infty}^{\infty}\lambda^nq_{\sigma B^{-n}}
\end{equation}
converges. The Hypothesis (H1) of \cite{kaplan1984} gives a sufficient condition to ensure this convergence.
\begin{definition}[Hypothesis H1]\label{H1}
We say that an almost periodic function $q\colon\mathbb{R}\to\mathbb{R}$ satisfies the hypothesis (H1) if
\begin{equation}
\tag{H1} \sum_{a} |q_a||a|^\alpha<\infty
\end{equation}
where $\alpha=-\log{\lambda}/\log{B}$ and $q_a$ are its Fourier Coefficients.
\end{definition}
Although this condition may appear to be too technical, \cite{kaplan1984} proved that if $p$ is a $C^3$ function then $q$ will satisfy (H1) (see Proposition \ref{propositionH1} for our version in higher dimensions).

\begin{lemma}
If $q$ is an almost periodic function satisfying $(\mathrm{H1})$, $0<\lambda<1$, and $B>1/\lambda$, then the series $(\ref{gsigma})$ converges.
\label{lemmaH1}
\end{lemma}
\begin{proof}
Given $|\sigma|\geq 1$ and $n\in\mathbb{N}$, let $a_n:=\sigma B^{-n}$, notice that $|a_n|^{\alpha}=|\sigma|^{\alpha}B^{-n\alpha}$, and since $\alpha=-\log{\lambda}/\log{B}$ we have
$B^{-n\alpha}=\lambda^n$, so
\[
|\lambda^nq_{\sigma B^{-n}}|=\left|\left(\frac{a_n}{\sigma}\right)^{\alpha}q_{a_n}\right|\leq |a_n|^\alpha |q_{a_n}|.
\]
Thus, Hypothesis (H1) ensures that 
\[
\sum_{n=-\infty}^{\infty}|\lambda^kq_{\sigma B^{-n}}|\leq \sum_{n=-\infty}^{\infty}  |a_n|^\alpha |q_{a_n}| \leq \sum_{a} |a|^\alpha |q_{a}|<\infty.
\]
If $|\sigma|<1$, then consider $k \in \mathbb{N}$ such that $|\sigma B^k|>1$ and note that
\begin{align*}
g_\sigma &=\sum_{n=-\infty}^{\infty}\lambda^nq_{\sigma B^{-n}}=\sum_{n=-\infty}^{\infty}\lambda^nq_{(\sigma B^k) B^{-n-k}} \\
&=\frac{1}{\lambda^k}\sum_{n=-\infty}^{\infty}\lambda^{n+k}q_{(\sigma B^k) B^{-n-k}} 
= \frac{1}{\lambda^k}\sum_{m=-\infty}^{\infty}\lambda^{m}q_{\tilde{\sigma}B^{-m}},  
\end{align*}
where $\tilde{\sigma}=\sigma B^k$, so we are in the previous case.
\end{proof}

So now that $g_\sigma$ is well defined, there are two possible cases: either (1) $g_\sigma=0$ for every $\sigma$, or (2) there is a Fourier coefficient $g_\sigma \neq 0$. In case (1), hypothesis (H1) ensures that $g(t)\equiv 0$ as shown in Proposition 2.2 from \cite{kaplan1984}, hence,
\[
f(t)=-\sum^{-1}_{n=-\infty} \lambda^n q(B^n t),
\]
which ensures that $f$ is $C^1$ by a similar argument as in Proposition \ref{dim1} (see Lemma 2.1 in \cite{kaplan1984}).
In case (2), $g_\sigma \neq 0$ can be used to ensure the existence of $C_1>0$ such that
\[
\underset{J}{\text{var}}(f) \geq C_1 L^{\alpha}
\]
for every interval $J$ with $L=\text{length}(J)$ small enough. Indeed, the following proposition is the main step in proving this.

\begin{proposition}[Proposition 2.5 in \cite{kaplan1984}]
    Given a continuous real function $\gamma(t)$, let 
    \[
    I=\frac{\rho}{2\pi n} \int_{0}^{\frac{2\pi n}{\rho}} \gamma(t)\cos(\rho t +\phi) dt.
    \]
    where $\rho\in(0,\infty)$, $\phi \in \mathbb{R}$ and $n$ is a positive integer. Then
    \[
    \underset{J}{\mathrm{var}}(\gamma) \geq \pi |I|,
    \]
    where $J=[0, 2\pi n/\rho]$.
    \label{prop2.5}
\end{proposition}

\citet{kaplan1984} use this proposition for any translation of the form $\gamma(t)=f(t+\theta)$, since the variation of $f(t+\theta)$ on $[0,2\pi n/\rho]$ equals the variation of $f(t)$ on $[\theta,2\pi n/\rho+\theta]$, so they can cover any interval of the domain. Using  Lemma 2.6 and Lemma 2.7 in \cite{kaplan1984} they shown that for large enough $k,n\in\mathbb{N}$ there is $0<C_0<|g_{\sigma_0}|$ such that
\[
|I|=\frac{\sigma_0B^k}{2\pi n}\int_0^{\frac{2\pi n}{\sigma_0B^k}} f(t+\theta)\cos(\sigma_0B^k(t+\theta)) dt \geq \frac{\lambda^k(|g_{\sigma_0}|-C_0)}{\pi}
\]
Setting $C_1$ as $|g_{\sigma_0}|-C_0$ and using Proposition \ref{prop2.5} we obtain
\begin{eqnarray*}
\underset{J}{\text{var}}(f)&\geq& \lambda^k C_1=B^{-\alpha k}C_1 \\
&\geq& \left(\frac{\text{length}(J)}{(2\pi n/\sigma_0)}\right)^{\alpha}C_1 \\
&=&C_1(\text{length}(J))^\alpha,
\end{eqnarray*}
whenever $\text{length}(J)\leq(2\pi n/\sigma_0)B^{-k}$.
For the upper bound, the fact that $q$ and $q'$ are bounded let us obtain $C_2>0$ such that
\[
\underset{J}{\text{var}}(f) \leq C_2 L^{\alpha}
\]
if $L=\text{length}(J)$ is small enough (see Proposition 2.9 in \cite{kaplan1984}). Then Proposition \ref{var} ensures that $$\mathrm{dim(graph}(f))=2-\alpha.$$
After calculating the box dimensions of the stable and unstable slices, it is possible to calculate the box dimension of the graph of $\phi$. We expand some details in the proof given in \cite{kaplan1984}.
\begin{theorem}[Theorem B in \cite{kaplan1984}]
Let $F$ be the skew product $($\ref{F}$)$ with $k=2$ and $\phi$ be the function defined in $($\ref{phi}$)$. If $p$ is $C^3$ and $\lambda \in (B_1,1)$, then either $\phi$ is nowhere differentiable and
\[
\mathrm{dim(graph}(\phi))=3-\left|\frac{\log \lambda}{\log B_1}\right|
\]
or $\phi$ is $C^1$ and $\mathrm{dim(graph}(\phi))=2$.
\label{kaplanthmB}
\end{theorem}

\begin{proof}
Recall that $\phi$ restricted to the unstable manifold is a $C^1$ function (see Proposition \ref{maior}) and when all Fourier coefficients of $g_1$ are null, $\phi_1$ is $C^1$, so in this case $\phi$ is $C^1$ and $\mathrm{dim(graph}(\phi))=2$. If there exists a non-zero Fourier coefficient of $g_1$, then there exist $C_1,C_2>0$, as in the proof of Theorem \ref{kaplanthmA}, and $K_1>0$ given by Proposition \ref{maior}, such that
\[
C_1L^{\alpha}\leq \underset{ |t|\leq \frac{L}{2}}{\text{var}}\phi_1\leq C_2L^\alpha \,\,\,\,\,\, \text{and} \,\,\,\,\,\, \underset{\mathbb{T}^2}{\sup}\left|\frac{\partial \phi}{\partial s}\right|<K_1,
\]
if $L$ is sufficiently small. If $(t_0,s_0)\in\mathbb{T}^2$ and $E_L:=\left[t_0-\frac{L}{2},t_0+\frac{L}{2}\right]\times\left[s_0-\frac{L}{2},s_0+\frac{L}{2}\right]$, then the above inequalities ensure that
$$C_1L^\alpha-2K_1L\leq  \underset{E_L}{\text{var}} \phi \leq C_2L^\alpha+2K_1L$$
that is equivalent to
$$L^{\alpha}(C_1-2K_1L^{1-\alpha})\leq \underset{E_L}{\text{var}} \phi \leq L^\alpha(C_2+2K_1L^{1-\alpha})$$ 
that, in turn, implies
$$\frac{1}{2}C_1L^\alpha \leq \underset{E_L}{\text{var}} \phi \leq 2C_2 L^\alpha$$
if $L^{1-\alpha}\leq \text{min}\{C_2/2K_1,C_1/4K_1\}$. Proposition \ref{nvar} ensures that
\[
\mathrm{dim(graph}(\phi))=3-\alpha = 3-\left|\frac{\log \lambda}{\log(1/B_1)}\right|= 3-\left|\frac{\log \lambda}{\log B_1}\right|.
\]
This completes the proof.
\end{proof}

\section{Fractal structure of the graph and its sub-slices}\label{Sec3}

In this section, we consider the case $k>2$, where $A$ is an Anosov automorphism with $k$ distinct eigenvalues $B_1,\dots,B_k$ satisfying $$|B_1|<|B_2|<\cdots<|B_j|<1<|B_{j+1}|<\cdots<|B_k|$$ and such that the respective eigenvectors $v_1,\dots,v_k$ form a basis of $\mathbb{R}^k$. Note that $j$ is the number of stable eigenvalues of $A$, that is, the dimension of the stable manifolds of $A$. In the whole section, the numbers $k,j,B_i$, vectors $v_i$, and the respective functions $q_i$, $\phi_i$, $g_i$ defined below in $($\ref{qi},\ref{phii},\ref{gi}$)$ will appear in statements of results repeatedly without proper explanation in the statement. We hope this causes no confusion.

\subsection{\texorpdfstring{Box dimension for the stable slices in $\mathbb{T}^k \times \mathbb{R}$}{Box dimension for the stable slices in Tk X R}}

\label{sub-slices}

In Proposition \ref{maior} we proved that if $\lambda<|B_i|$ then the box dimension of the graph of $\phi_i$ is 1. When $|B_i|<\lambda$, the function $\phi_i$ can be written in the form of the function $f$ in Theorem $\ref{kaplanthmA}$. Indeed, letting 
\begin{equation}\label{qi}
q_i(t)=p(tv_i+\underline{t})
\end{equation}
it follows that 
\begin{equation}\label{phii}
\phi_i(t)=\sum_{n=0}^{\infty}\lambda^n q_i((B_i^{-1})^{n} t).
\end{equation}
Thus, the argument of subsection 2.4 can also be applied for $\phi_i$, using the function 
\begin{equation}\label{gi}
g_i(t):=\sum_{n=-\infty}^{\infty} \lambda^{n}q_i((B_i^{-1})^{n} t),
\end{equation}
associated to $q_i$ as in (\ref{g}), and we obtain the following theorem.

\begin{theorem}
Let $F$ be the skew product $($\ref{F}$)$ and for each $i\leq j$, let $q_i$, $\phi_i$, $g_i$ be the functions defined in $($\ref{qi},\ref{phii},\ref{gi}$)$.
If $q_i$ is a $C^1$ function, almost periodic, and satisfies the hypothesis $(\mathrm{H1})$, then either all the Fourier coefficients of $g_i$ are zero, $\phi_i$ is $C^1$, and $\mathrm{dim(graph}(\phi_i))=1$, or 
there is at least one non-zero Fourier coefficient of $g_i$ and
\[
\mathrm{dim(graph}(\phi_i)) = 2-\frac{\ln{\lambda}}{\ln{|B_i|}}.
\]
\end{theorem}
The hypothesis on $q_i$ are satisfied with enough regularity on $p$. Indeed, $q$ will be almost periodic as long as $p$ is continuous and the Hypothesis (H1) is satisfied when $p$ is $C^{k+1}$. This is proved in the next results but analogous results in the case $k=2$ are announced in \cite{kaplan1984}.

\begin{lemma}
    If $p$ is uniformly continuous, then $q_i$ is almost periodic.  
\end{lemma}
\begin{proof}
Given $\epsilon>0$, because $p$ is uniformly continuous, there is $\delta>0$ such that 
\[
d(x,y)\leq \delta \Rightarrow |p(x)-p(y)|<\epsilon.
\]
Let $\alpha_i(t)=\underline{t}+tv_i$ be a parametrization of the stable manifold on the direction of $v_i$. Since $A$ is linear, $\alpha_i$ parametrizes an orbit of a linear flow on the Torus $\mathbb{T}^k$, that is a recurrent flow without singularities. Thus, if $\alpha_i(t)\in \overline{B(\underline{t},\delta)}$, we can consider the smaller $s>0$ satisfying: $\alpha_i(t+s)\in\overline{B(\underline{t},\delta)}$ and such that there exists $s'\in(t,s)$ with $\alpha_i(t+s')\notin\overline{B(\underline{t},\delta)}$.
We let $s=\tau(\alpha_i(t))$ and call it the first-return time of $\alpha_i(t)$ to $\overline{B(\underline{t},\delta)}$.
Let 
\[
L=\sup\{\tau(\alpha_i(t)); \,\,\alpha_i(t)\in\overline{B(\underline{t},\delta)}\} 
\] and note that $L<\infty$.
Every interval $I$ of length $L$ contains $s\in I$ satisfying
\[
d(\alpha_i(0),\alpha_i(s))\leq \delta,  
\]
since $d(\alpha_i(0),\alpha_i(0+s))>\delta$ for every $s\in I$ implies the existence of $t\in\mathbb{R}$ such that $$\alpha_i(t)\in\overline{B(\underline{t},\delta)} \,\,\,\,\,\, \text{and} \,\,\,\,\,\, \tau(\alpha_i(t))>L$$ contradicting the definition of $L$.
Now given $r,s,t\in \mathbb{R}$, by linearity we have
\[
d(\alpha_i(t),\alpha_i(t+s))=d(\alpha_i(r),\alpha_i(r+s)),
\]
and this implies that for each interval $I$ of length $L$, there is $s\in I$ such that 
\[
d(\alpha_i(t),\alpha_i(t+s))\leq \delta \,\,\,\,\,\, \text{for every} \,\,\,\,\,\, t\in\mathbb{R}.
\]
Hence,
\[
|q_i(t)-q_i(t+s)|=|p(\alpha_i(t))-p(\alpha_i(t+s))|<\epsilon.
\]
for every $t\in\mathbb{R}$. Since this can be done for each $\epsilon>0$, the proof is complete.
\end{proof}

\begin{proposition}
\label{propositionH1}
If $p\colon\mathbb{T}^k\rightarrow\mathbb{R}$ is $C^{k+1}$, then $q_i$ is almost periodic and satisfies $(\mathrm{H1})$ for every $i\in\{1,\dots,j\}$.
\label{pck1}
\end{proposition}
\begin{proof} 
 As $p$ is periodic and $C^{k+1}$, the Fourier Series of $p$ is convergent, let
 \[
 p(\underline{x})=\sum_{\underline{j}\in \mathbb{Z}^k}p_{\underline{j}}e^{2\pi i(\underline{x}\cdot\underline{j})}.
 \]
 Because $p$ is $C^{k+1}$, we have that there is a positive $K$ such that
 \[
 |p_{\underline{j}}|\leq \frac{K}{|\underline{j}|^{k+1}}.
 \]
Indeed, the coefficient of $\partial^{k+1} p/\partial i^{k+1}$ are $2\pi (j_i)^{k+1}p_{\underline{j}}$ and they are all limited (see \cite{corduneanu}) therefore 
\[
|p_{\underline{j}}|\leq \frac{1}{k}\sum_{i=1}^k \frac{K_i}{2\pi |j_i|^{k+1}}\leq \frac{ \sum_{i=1}^k K_i}{2\pi k\prod_{i=1}^k |j_i|^{k+1}} 
\]
because $|j_i|\geq 1$ we have 
\[
 k\prod_{i=1}^k |j_i|\geq \sum_{i=1}^k \sqrt{|j_i|^2} 
\]
and so 
\[
|p_{\underline{j}}|\leq \frac{ \sum_{i=1}^k K_i}{2\pi \sum_{i=1}^k \sqrt{ |j_i|^2} } \leq \frac{K}{|\underline{j}|^{k+1}}.
\]
With this the Fourier series of $q_i(t)$ can be written as
 \begin{equation}
 \sum_{\underline{j}\in \mathbb{Z}^k}p_{\underline{j}}e^{2\pi i((tv_i+\underline{t})\cdot\underline{j})}=\sum_{\underline{j}\in \mathbb{Z}^k}e^{2\pi i(\underline{t}\cdot\underline{j})} p_{\underline{j}}e^{2\pi i((tv_i)\cdot\underline{j})}=\sum_{a}(q_{i})_{a} e^{iat}
\label{relacaopeq}
\end{equation}
where $a=2\pi(\underline{j}\cdot v_i)$.
Thus, using the Cauchy-Schwarz inequality and the one above,
\begin{align*}
    \sum_{a}|(q_{i})_{a}||a|^{\alpha}&=    (2\pi)^{\alpha}\sum_{\underline{j}\in \mathbb{Z}^k} |e^{2\pi i \underline{t}\cdot\underline{j}}p_{\underline{j}}||\underline{j}\cdot v_i|^\alpha= 
   (2\pi)^{\alpha}\sum_{\underline{j}\in \mathbb{Z}^k}|p_{\underline{j}}||\underline{j}\cdot v_i|^\alpha \\ 
   &\leq (2\pi)^{\alpha}\sum_{\underline{j}\in \mathbb{Z}^k}\frac{K|\underline{j}|^\alpha|v_i|^\alpha}{|\underline{j}|^{k+1}} 
     \leq K_2 \sum_{\underline{j}\in \mathbb{Z}^k} \frac{1}{|j|^{k+1-\alpha}}.
\end{align*}
This converges because $k+1-\alpha>k$ (see chapter 5, section 32 of \cite{bromwich2005}\footnote{This fact is proved for the double sum, but also works for any finite number sum.}). 
\end{proof}

The following is a direct corollary of the above results.

\begin{theorem}
Let $F$ be the skew product $($\ref{F}$)$ and for each $i\leq j$, let $\phi_i$, $g_i$ be the functions defined in $($\ref{phii},\ref{gi}$)$. If $p$ is a $C^{k+1}$ function and $i\leq j$, then either all the Fourier coefficients of $g_i$ are zero, $\phi_i$ is $C^1$, and $\mathrm{dim(graph}(\phi_i))=1$, or 
there is at least one non-zero Fourier coefficient of $g_i$ and
\[
\mathrm{dim(graph}(\phi_i)) = 2-\frac{\ln{\lambda}}{\ln{|B_i|}}.
\]
\label{thm1}
\end{theorem}

This calculates the box dimension of all stable sub-slices of the graph.

\subsection{Dimension of the graph and dimension of the sub-slices}

In this subsection, we discuss how to calculate the dimension of the graph using the dimension of the sub-slices calculated in the previous section and compare the dimension of the graph with the sum of the dimensions of its sub-slices. We obtain sufficient conditions that ensure the dimension of the graph is smaller than the sum of the dimension of the sub-slices (see Theorem \ref{thm2}). This means that there exist fractal structure in each sub-slice that is not captured by the box dimension of the graph. The box dimension of the graph in the case $k>2$ can be calculated as in Theorem C of \cite{kaplan1984}. Indeed, if $|B_1|<\lambda<1$ and there is $\sigma$ such that $(g_1)_{\sigma}\neq 0$, then the variation of $\phi$ is proportional to the variation of $\phi_1$ and so
\[
\mathrm{dim(graph}(\phi))=k+1-\frac{\log \lambda}{\log |B_{1}|}=k-1+\mathrm{dim}(\phi_1).
\] 
The above number is also the Lyapunov dimension of the $\mathrm{graph}(\phi)$ (see \cite{kaplan83} for more details), and it is conjectured that this is usually the case for an attractor.\footnote{You can find more results in this direction in \cite{Ledrappier}.} If $(g_1)_{\sigma}=0$ for every $\sigma$ and $|B_2|<\lambda<1$, we can try to see if there is fractal dimension for the sub-slice $\phi_2$ by seeing if there is $\sigma$ such $(g_2)_{\sigma}\neq 0$. If this is the case, then 
\[
\mathrm{dim(graph}(\phi))=k+1-\frac{\log \lambda}{\log |B_{2}|}=k-1+\mathrm{dim}(\phi_2).
\] 
and if not, we can repeat the argument till the sub-slice $\phi_l$ as long as $|B_l|<\lambda$. The idea is that the variation of $\phi$ and the variation of the sub-slice with the highest variation are proportional to the same value, while the sub-slice with the highest variation is given by the smallest $i$ admitting a $(g_i)_{\sigma}\neq0$.

\begin{theorem}
Let $F$ be the skew product $($\ref{F}$)$ and $\phi$ be the function defined in $($\ref{phi}$)$.
If $p$ is $C^{k+1}$, $\lambda \in (B_{l},B_{l+1})\cap (0,1)$ with $l\leq j$, and there is $(g_i)_\sigma \neq 0$ with $i\leq l$, then
\[
\mathrm{dim(graph}(\phi))=k+1-\frac{\log \lambda}{\log |B_{i_0}|},
\]
where $i_0=\mathrm{min}\{i\leq l \, \, | \, \, \exists\,\, \sigma, (g_{i})_\sigma \neq 0 \}$. If $(g_i)_{\sigma}=0$ for every $\sigma$ and every $i\in\{1,\dots,j\}$, then $\phi$ is $C^1$ and $\mathrm{dim(graph}(\phi))=k$.
\label{kaplanthmC}
\end{theorem}

\begin{proof}
For each $i\in\{1,\dots,l\}$, let
\[
\alpha_i=\begin{cases}
    \frac{\log \lambda}{\log|B_i|}, & \text{if there is $\sigma$ such that } (g_{i})_{\sigma}\neq 0\\
    1, & \text{otherwise}
\end{cases}
\]
and if $i>l$ let $\alpha_i=1$. Consider the $k-$cube $I=\left[-\frac{L}{2},\frac{L}{2}\right]^k$, so we have that for each $i\in\{1,\dots,j\}$ there are positive constants $C_i,K_i$ such that
\[
K_iL^{\alpha_i}\leq \underset{\left[-\frac{L}{2},\frac{L}{2}\right]}{\text{var}}(\phi_i) \leq C_iL^{\alpha_i}.
\]
Since 
$$\underset{I}{\text{var}}(\phi)\geq\underset{\left[-\frac{L}{2},\frac{L}{2}\right]}{\text{var}}(\phi_i) \,\,\,\,\,\, \text{for every} \,\,\,\,\,\, i\in\{1,\dots,k\}$$ and $$\underset{I}{\text{var}}(\phi)\leq\sum_{i=1}^k\underset{\left[-\frac{L}{2},\frac{L}{2}\right]}{\text{var}}(\phi_i),$$ it follows that
\[
K_{i_0}L^{\alpha_{i_0}}\leq \underset{I}{\text{var}}(\phi) \leq \sum^{k}_{i=1}C_iL^{\alpha_i}.
\]
Note that $\alpha_i-\alpha_{i_0}>0$ for every $i\neq i_0$. Indeed, this is clear if $\alpha_i=1$, and otherwise, $i_0< i\leq l$ and 
$$\alpha_i=\frac{\log \lambda}{\log|B_i|}>\frac{\log \lambda}{\log|B_{i_0}|}=\alpha_{i_0}.$$
Thus, if $L$ is small enough such that
\[
L^{\alpha_i-\alpha_{i_0}}\leq\frac{C_{i_0}}{C_i} \,\,\,\,\,\, \text{for every} \,\,\,\,\,\, i\neq i_0,
\]
then
\[
K_{i_0}L^{\alpha_{i_0}}\leq \underset{I}{\text{var}} (\phi) \leq kC_{i_0}L^{\alpha_{i_0}}.
\]
So by the Proposition \ref{nvar} we have 
\[
\mathrm{dim(graph}(\phi))=k+1-\frac{\log \lambda}{\log |B_{i_0}|}=k-1+\mathrm{dim(graph}(\phi_{i_0})).
\]
\end{proof}

It is interesting that the only fractal structure captured by the box dimension of the graph is the one given by the sub-slice $\phi_{i_0}$. Thus, if there are more fractal sub-slices $\phi_i$ for $i>i_0$, then the sum of the box dimensions of the sub-slices will be greater than the box dimension of the graph. The following is the main result of this article.

\begin{theorem}\label{thm2}
Let $F$ be the skew product $($\ref{F}$)$ and $\phi$ be the function defined in $($\ref{phi}$)$.
If $p$ is $C^{k+1}$, $\lambda \in (B_{l},B_{l+1})\cap (0,1)$ with $l\leq j$, and there exist $i_1\in\{i_0+1,\dots,l\}$ and $\sigma_1\in\mathbb{R}$ with $(g_{i_1})_{\sigma_1}\neq0$, then 
\[
\mathrm{dim(graph}(\phi))<\sum_{i=1}^{k}\mathrm{dim(graph}(\phi_i)).
\]
\end{theorem}

\begin{proof}
Theorem \ref{kaplanthmC} ensures that
\[
\mathrm{dim(graph}(\phi))=k-1+\mathrm{dim(graph}(\phi_{i_0})),
\]
Proposition \ref{dim1} ensures that 
$$\mathrm{dim(graph}(\phi_m))=1 \,\,\,\,\,\, \text{for every} \,\,\,\,\,\, m\in\{l+1,\dots,k\},$$ 
and by the definition of $i_0$, Theorem \ref{thm1} ensures that the same is true for $m<i_0$.
The hypothesis and Theorem \ref{thm1} ensure that
\[
\mathrm{dim(graph}(\phi_{i_1}))=2-\frac{\log\lambda}{\log |B_{i_1}|} > 1.
\]
Since
\[
\mathrm{dim(graph}(\phi_m))\geq 1 \,\,\,\,\,\, \text{whenever} \,\,\,\,\,\, i_0<m\leq l,
\]
it follows that
\[
\sum_{m\neq i_0} \mathrm{dim(graph}(\phi_m))>k-1
\]
and, hence,
\[
\mathrm{dim(graph}(\phi))=k-1+\mathrm{dim(graph}(\phi_{i_0}))<\sum_{i=1}^{k}\mathrm{dim(graph}(\phi_l)).
\]
\end{proof}

\subsection{Examples and generic conditions}

In this subsection, we prove that the hypothesis on Theorem \ref{thm2} about the existence of two distinct $i_0,i_1\in\{1,\dots,l\}$ with $(g_{i_0})_{\sigma_0}\neq0$ and $(g_{i_1})_{\sigma_1}\neq0$ is always satisfied if we choose correctly the function $p\colon\mathbb{T}^k\to\mathbb{R}$. Actually, we prove that it is possible to choose $p$ such that this holds for every $i\in\{1,\dots,l\}$ and that the set of $C^{k+1}$ functions satisfying this is generic in the $C^{k+1}$ topology. We  begin with an explicit example illustrating Theorem \ref{thm2}. 
\begin{example}\label{example}
Let A be the following hyperbolic matrix
\begin{align*}
\begin{pmatrix}
6 & -5 & 1 \\
1 & 0 & 0 \\
0 & 1 & 0 
\end{pmatrix}.
\end{align*}
The eigenvalues $\mu_1,\mu_2,\mu_3$ of $A$ satisfy $$0<\mu_1<\mu_2<0.65<1<5<\mu_3.$$

Moreover, if $\mu$ is an eigenvalue of $A$, then $v_{\mu}= (\mu^2, \mu, 1)$ is an eigenvector associated to $\mu$. For each $i\in\{1,2,3\}$ let $v_i=(\mu_i^2,\mu_i,1)$. 
For each $i\in\{1,2\}$, let $a_i>0$, $\theta_{i} \in (0,\pi/2)$ and
\begin{equation}
q_{i}(t)=a_{i} \cos{(\mu_i^{-1}t+\theta_{i})}.
\label{funcaoQ}    
\end{equation}

\begin{proposition}
In the coordinate system given by $(v_1,v_2,v_3)$, if
\[
p(x_1,x_2,x_3):=q_{1}(x_1) +q_{2}(x_2) +  q_3(x_3),
\]
where $q_3$ is any $C^4$ function, then for each $i\in\{1,2\}$ we have
\[
(g_i)_{1}=\lambda(q_i)_{\mu_i^{-1}}\neq 0.
\]
\label{p}
\end{proposition}
\begin{proof}
This stands from the fact that the Fourier coefficient $(q_i)_{\mu_i^{-1}}$ is the only non zero coefficient of $q_{i}$. Indeed, it is known that
\[
\lim_{T\rightarrow \infty} \frac{1}{2T}\int_{-T}^{T} \cos{(bt+\theta)}\sin(at)\hspace{0.5em} dt=\begin{cases} 0 & \text{ if } b\neq a \\
\frac{-\sin(\theta)}{2}& \text{ if } b= a\end{cases}
\]
and
\[
\lim_{T\rightarrow \infty} \frac{1}{2T}\int_{-T}^{T} \cos{(bt+\theta)}\cos(at)\hspace{0.5em} dt=\begin{cases} 0 & \text{ if } b\neq a \\
\frac{\cos(\theta)}{2}& \text{ if } b= a.\end{cases}
\]
Therefore,
\begin{align*}
(q_i)_a=\lim_{T\rightarrow\infty}\frac{1}{2T}\int^{T}_{-T} a_i\cos{(\mu_i^{-1}t+\theta_{i})} e^{-iat}\hspace{0.5em} dt
\end{align*}
is non zero if, and only if, $a=\mu_i^{-1}$, and in this case, $(q_i)_{\mu_i^{-1}}=\frac{a_ie^{i\theta_i}}{2}$. By Definition (\ref{gsigma}) it follows that 
\[
(g_i)_{1}=\lambda(q_i)_{\mu_i^{-1}}\neq 0.
\]
\end{proof}
Thus, if $\lambda \in (\mu_2,1)$, then 
\[
\mathrm{dim(graph}(\phi_1))=2-\frac{\log\lambda}{\log \mu_1}, \,\,\,\,\,\,   \mathrm{dim(graph}(\phi_2))=2-\frac{\log\lambda}{\log \mu_2},
\]
and, hence,
\[
\mathrm{dim(graph}(\phi))=4-\frac{\log\lambda}{\log \mu_1}<\sum_{i=1}^{3}\mathrm{dim(graph}(\phi_i))=\left(2-\frac{\log\lambda}{\log \mu_1}\right)+\left(2-\frac{\log\lambda}{\log \mu_2}\right)+1.
\]
If $\lambda<\mu_2$, then
\[
\mathrm{dim(graph}(\phi))=\sum_{i=1}^{3}\mathrm{dim(graph}(\phi_i)).
\]
Indeed, if $\lambda \in (\mu_1,\mu_2)$, then
\[
\mathrm{dim(graph}(\phi_1))=2-\frac{\log\lambda}{\log \mu_1}, \,\,\,\,\,\,   \mathrm{dim(graph}(\phi_2))=1,
\]
and, hence,
\[
\sum_{i=1}^{3}\mathrm{dim(graph}(\phi_i))=\left(2-\frac{\log\lambda}{\log \mu_1}\right)+1+1=4-\frac{\log\lambda}{\log \mu_1}=\mathrm{dim(graph}(\phi)),
\]
and if $0<\lambda<\mu_1$, then $\phi$ is $C^1$ and $\mathrm{dim(graph}(\phi))=3$.
\end{example}
Now we return to the general case where $A$ is any Anosov automorphism of the Torus as in the beginning of this section.
\begin{theorem}\label{3.9}
Let $F$ be the skew product $($\ref{F}$)$ and $\phi$ be the function defined in $($\ref{phi}$)$.
If $\lambda \in (|B_l|, |B_{l+1}|)\cap(0,1)$ with $l\geq 2$, then there is a $C^{k+1}$ function $p$ satisfying: for each $i\in \{1,\dots,l\}$ there exists $\sigma_i\in\mathbb{R}$ such that $(g_{i})_{\sigma_i}\neq 0$. Moreover, the set of functions satisfying this is an open and dense subset in $C^{k+1}$ (in the $C^{k+1}$ topology).
\end{theorem}
\begin{proof}
The first part of this Theorem follows the idea of Proposition \ref{p}.
If $p$ is any $C^{k+1}$ function, let $I$ be the set of $i\in\{1,\dots,l\}$ such that $(g_{i})_{\sigma}=0$ for every $\sigma\in\mathbb{R}$. For each $\epsilon>0$, the function
\[
p(\underline{x})+\epsilon\left(\sum_{i\in I}\cos(\mu_i^{-1}x_i)\right)
\]
satisfies $(g_{i})_1\neq 0$ for every $i\in I$ (by Proposition \ref{p}). Also, if $i\notin I$, it satisfies $(g_{i})_{\sigma_i}\neq0$ for some $\sigma_i\in\mathbb{R}$. Letting $\epsilon\to0$ we obtain that $p$ is approximated by $C^{k+1}$ functions satisfying: for each $i\in \{1,\dots,l\}$ there exists $\sigma_i\in\mathbb{R}$ such that $(g_{i})_{\sigma_i}\neq 0$. Because of the relation (\ref{relacaopeq}) we have that the Fourier coefficients of $q_i$ depend continuously on the Fourier coefficients of $p$ and so $(g_i)_\sigma$ also varies continuously with $p$. This concludes the proof.
\end{proof}
We conclude noting that for a generic function $p$ the box dimension of a stable sub-slice $\phi_i$ will depend only if $\lambda<|B_i|$ or $\lambda>|B_i|$. In the first case, the dimension is one and in the second is $2-\log{\lambda}/\log{|B_i|}$.  
\begin{corollary}
Let $F$ be the skew product $($\ref{F}$)$ and $\phi$ be the function defined in $($\ref{phi}$)$. If $p$ is a generic $C^{k+1}$ function and $\lambda \in (|B_l|, |B_{l+1}|)\cap(0,1)$ for some $l\leq j$, then
    \[
\sum_{i=1}^{k}\mathrm{dim(graph}(\phi_i))=(k-l)+2l-\sum_{i=0}^{l} \frac{\log{\lambda}}{\log{|B_i|}}.
    \]
\end{corollary} 

\section*{Acknowledgments} The authors are in debt with the anonymous referee for careful reading of this article and for suggestions that considerably improved its presentation.

\section*{Funding}
Bernardo Carvalho was supported by Progetto di Eccellenza MatMod@TOV grant number PRIN 2017S35EHN and by CNPq grant numbers 405916/2018-3 and 446192/2024. Both authors were supported by Fapemig grant number APQ-00036-22.
\nocite{Keller}
\nocite{Takens}
\nocite{Palis-Viana}
\nocite{BARANSKI201432}
\nocite{Shen}
\nocite{HADJILOUCAS_NICOL_WALKDEN_2002}
\nocite{Walkden}
\bibliographystyle{abbrvnat}
\bibliography{biblioteca}
\end{document}